\documentclass[11pt]{amsart}

\usepackage{amssymb}
\input epsf

\newtheorem{thm}{Theorem}[section]
\newtheorem{lemma}[thm]{Lemma}

\def\im{\mathop{\rm Im}\nolimits \,}

\newcommand{\cL}{{\mathcal L}}

\begin{document}

\thanks{Partially supported by NSF grant DMS-0100108}

\title{On the Grinberg -- Kazhdan formal arc theorem}

\author{Vladimir Drinfeld}

\address{Dept. of Math., Univ. of Chicago, 5734 University 
Ave., Chicago, IL 60637}

\email{drinfeld@math.uchicago.edu}

\maketitle


Let $X$ be a scheme of finite type over a field $k$, and 
$X^\circ \subset X$ the smooth part of $X$.  Consider the scheme $\cL(X)$
of formal arcs in $X$.  The $k$-points of $\cL(X)$ are just maps
${\rm Spec} \, k[[t]] \to X$. Let $\cL^\circ(X)$ be the
open subscheme of arcs whose image is not contained in
$X \setminus X^\circ$.  Fix an arc $\gamma_0:{\rm Spec}\, k[[t]]\to X$ in
$\cL^\circ(X)$, and let $\cL(X)_{\gamma_0}$ be the formal
neighborhood of $\gamma_0$ in $\cL(X)$. We will give a simple proof of
the following theorem, which was proved by M.~Grinberg and D.~Kazhdan
\cite{GK} for fields $k$ of characteristic $0$.

\begin{thm}\label{main}
There exists a scheme $Y = Y(\gamma_0)$ of finite type over $k$,
and a point $y \in Y(k)$, such that
$\cL(X)_{\gamma_0}$ is isomorphic to $D^{\infty}\times Y_y$
where $Y_y$ is the formal neighborhood of $y$ in $Y$ and 
$D^{\infty}$ is the product of countably many copies of the formal disk
$D:={\rm Spf}\, k[[t]]$.
\end{thm}

\medskip

\noindent
{\bf A convention.}  Throughout this paper, a {\it test-ring} $A$ is a
local commutative unital $k$-algebra with residue field $k$ whose
maximal ideal $m$ is nilpotent. If $S$ is a scheme over $k$ and 
$p\in S(k)$ is a $k$-point we think of the formal neighborhood $S_p$ in
terms of its functor of points $A \mapsto S_p(A)$ from test-rings to
sets. For instance, $A$-points of $\cL(X)_{\gamma_0}$ are $A[[t]]$-points
of $X$ whose reduction modulo $m$ equals $\gamma_0$.

\medskip

\noindent
{\bf Remark.} 
$\cL(X)_{\gamma_0}=\cL(X_1)_{\gamma_0}$ where $X_1\subset X$ is the
closure of the connected component of $X^\circ$ containing 
$\gamma_0({\rm Spec} \, k((t))\,)$. So we can assume that $X$ is reduced
and irreducible. But $Y$ is, in general, neither reduced nor irreducible
(e.g., see the following example).

\medskip

\noindent
{\bf Example.} 
$X$ is the hypersurface $yx_{n+1}+f(x_1,\ldots,x_n)=0$ and $\gamma_0(t)$
is defined by $x_{n+1}^0(t)=t,y^0(t)=x_1^0(t)=\ldots=x_n^0(t)=0$. Then
one can define $Y$ to be the hypersurface $f(x_1,\ldots,x_n)=0$ and $y$
to be the point $0\in Y$. Indeed, by the Weierstrass division theorem
(alias preparatory lemma) for any test-ring $A$ every $A$-deformation of 
$x_{n+1}^0(t)=t$ has the form $x_{n+1}(t)=(t-\alpha)u(t)$ where $\alpha$
belongs to the maximal ideal $m\subset A$ and $u\in A[[t]]$ is
invertible. Given $\alpha$, $u$, and $x_1(t)\ldots,x_n(t)\in m[[t]]$
there is at most one $y(t)\in m[[t]]$ such that
$y(t)x_{n+1}(t)+f(x_1(t),\ldots,x_n(t))=0$, and $y(t)$ exists if and
only if $f(x_1(\alpha ),\ldots,x_n(\alpha ))=0$.

\medskip

\noindent
{\bf Proof of Theorem \ref{main}.} 
We can assume that $X$ is a closed subvariety of an affine space. Then
there is a closed subscheme $X'$ of the affine space such that
$X'\supset X$, $X'$ is a complete intersection, and the image of
our arc $\gamma_0$ is not contained in the closure of $X'\setminus X$.
Clearly $\cL(X)_{\gamma_0}=\cL(X')_{\gamma_0}$, so we can assume that
$X=X'$ is the subscheme of 
${\rm Spec}\, k[x_1,\ldots,x_n,y_1,\ldots,y_l]$ 
defined by equations $p_1=\ldots =p_l=0$ such that the arc $\gamma_0(t)
=(x^0(t),y^0(t))=(x^0_1(t),\ldots,x^0_n(t),y^0_1(t),\ldots,y^0_l(t))$ is
not contained in the subvariety of $X$ defined by 
$\det\frac{\partial p}{\partial y}=0$. Here
$\frac{\partial p}{\partial y}$ is the matrix of partial derivatives
$\frac{\partial p_i}{\partial y_j}$. Let $\gamma$ be an $A$-deformation
of $\gamma _0$ for some test-ring $A$, so $\gamma (t)=(x(t),y(t))$, where
$x(t)\in A[[t]]^n$, $y(t)\in A[[t]]^l$. Then by the Weierstrass division
theorem (alias preparatory lemma) 
$\det\frac{\partial p}{\partial y}(x(t),y(t))$ has a unique
representation as $q(t)u(t)$ where $u\in A[[t]]$ is invertible and 
$q$ is a monic polynomial whose reduction modulo the maximal ideal
$m\subset A$ is a power of $t$. Let $d$ denote the degree of $q$; it
depends only on $\gamma _0$, not on its deformation $\gamma$. We assume
that $d>0$ (otherwise we can eliminate $y$). The idea is to consider $q$
as one of the unknowns. More precisely, $A$-deformations of $\gamma_0$ 
are identified with solutions of the following system of equations. The
unknowns are $q(t)\in A[t]$, $x(t)\in A[[t]]^n$, and $y(t)\in A[[t]]^l$
such that $q$ is monic of degree $d$, $q(t)$ is congruent to $t^d$ modulo
$m$, and the reduction of $(x(t),y(t))$ modulo $m$ equals
$\gamma_0(t)=(x^0(t),y^0(t))$. The equations are as
follows:
\begin{equation}\label{1}
  \det\frac{\partial p}{\partial y}(x(t),y(t))\equiv 0 \mbox{ mod }q, 
\end{equation}
\begin{equation}\label{2}
   p(x(t),y(t))=0,
\end{equation}
where $p:=(p_1,\ldots,p_l)$. (Notice that if (\ref{1}) is satisfied then 
$q(t)^{-1}\det\frac{\partial p}{\partial y}(x(t),y(t))$ is
automatically invertible because it is invertible modulo $m$).

Now fix $r\ge 1$ and consider the following system of equations.
The unknowns are $q(t)\in A[t]$, $x(t)\in A[[t]]^n$, and 
$\bar y\in A[t]^l/(q^r)$ such that $q$ is monic of degree $d$, $q(t)$ is
congruent to $t^d$ modulo $m$, the reduction of $x(t)$ modulo $m$ equals
$x^0(t)$, and the reduction of $\bar y$ modulo $m$ equals the reduction
of $y^0$ modulo $t^{r}$. The equations are as follows:
\begin{equation}\label{3}
\det\frac{\partial p}{\partial y}(x(t),\bar y)\equiv 0\mbox{ mod }q, 
\end{equation}
\begin{equation}\label{4}
p(x(t),\bar y)\in\im (q^r\frac{\partial p}{\partial y}(x(t),\bar y):
A[t]^l/qA[t]^l\to q^rA[t]^l/q^{r+1}A[t]^l).
\end{equation}
Condition (\ref{4}) makes sense because $p(x(t),\bar y)$ is well defined
modulo $\im q^r\frac{\partial p}{\partial y}(x(t),\bar y)$. Notice that
(\ref{4})~is indeed an equation because it is equivalent to the condition
$\hat Cp(x(t),y(t))\equiv 0\mbox{ mod }q^{r+1}$ where $y(t)\in A[t]^l$ is
a preimage of $\bar y$ and $\hat C$ is the matrix adjugate to
$C:=\frac{\partial p}{\partial y}(x(t),y(t))$ 
(so $C\hat C=\hat CC=\det C$). This condition is equivalent to the
following equations, which do not involve a choice of $y(t)\in A[t]^l$
such that $y(t)\mapsto \bar y$:
\begin{equation}\label{5}
p(x(t),\bar y)\equiv 0\mbox{ mod }q^r,
\end{equation}
\begin{equation}\label{6}
\hat B p(x(t) ,\bar y)\equiv 0\mbox{ mod }q^{r+1},
\end{equation}
where $B:=\frac{\partial p}{\partial y}(x(t),\bar y)$; notice that
(\ref{6}) makes sense as soon as (\ref{5}) holds.

\begin{lemma}
The natural map from the set of solutions of (\ref{1}-\ref{2}) to the
set of solutions of (\ref{3}-\ref{4}) is bijective.
\end{lemma}

\noindent
{\bf Proof.} Let $a$ be the minimal number such that $m^a=0$. We proceed
by induction on $a$, so we can assume that $a\ge 2$ and the lemma is
proved for $A/m^{a-1}$. Then there exists $\tilde y(t)\in A[t]^l$ such
that $\tilde y\mbox{ mod }q^r=\bar y$ and 
$p(x(t),\tilde y(t))\in m^{a-1}[t]^l$; such $\tilde y$ is unique modulo
$q^rA[t]^l\cap m^{a-1}[t]^l$. We have to find 
$z(t)\in q^rA[t]^l\cap m^{a-1}[t]^l$ such that
$p(x(t),\tilde y(t)-z(t))=0$, i.e., $Cz(t)=p(x(t),\tilde y(t))$, where 
$C:=\frac{\partial p}{\partial y}(x(t),\tilde y(t))$. (\ref{3}) implies
that $\det C=q(t)u(t)$ for some invertible $u\in A[t]$. So $z(t)$ is
unique. By (\ref{4}) 
$p(x(t),\tilde y(t))\in q^rCA[t]^l+q^{r+1}A[t]^l$. But 
$CA[t]^l\supset (\det C)A[t]^l=qA[t]^l$, so $p(x(t),\tilde y(t))=Cz(t)$
for some $z(t)\in q^rA[t]^l$. We have 
$Cz(t)=p(x(t),\tilde y(t))\equiv 0\mbox{ mod }m^{a-1}$, so 
$q(t)z(t)\equiv 0\mbox{ mod }m^{a-1}$ and finally
$q(t)\equiv 0\mbox{ mod }m^{a-1}$.
\hfill\qedsymbol

\medskip

So the set of $A$-deformations of $\gamma_0$ can be identified with the
set of solutions of the system (\ref{3}-\ref{4}). This system is
essentially finite because $x(t)$ is relevant only modulo $q^{r+1}$.
E.g., if $r=1$ we can write
$x(t)$ as $q^2(t)\xi (t)+\bar x$, $\xi\in A[[t]]^n$, $\bar x\in A[t]^n$,
$\deg\bar x<2d$, and consider $\xi (t)$, $\bar x$, $q(t)$, and $\bar y$
to be the unknowns (rather than $x(t)$, $q(t)$, $\bar y$); then
(\ref{3}-\ref{4}) becomes a finite system of equations for $q$, $\bar x$,
$\bar y$ (and $\xi$ is not involved in these equations). So
$\cL(X)_{\gamma_0}$ is isomorphic to $D^{\infty}\times Y_y$, where the
$k$-scheme $Y$ of finite type and the point $y\in Y(k)$ are defined as
follows: for every $k$-algebra $R$ the set $Y(R)$ consists of triples 
$(q,\bar x ,\bar y)$ where $q\in R[t]$ is monic of degree $d$, 
$\bar x\in R[t]^n/(q^2)$, $\bar y\in R[t]^l/(q)$,
$\det B\equiv 0\mbox{ mod }q$, 
$B:=\frac{\partial p}{\partial y}(\bar x,\bar y)$, 
$p(\bar x,\bar y)\equiv 0\mbox{ mod }q$, and 
$\hat B p(\bar x ,\bar y)\equiv 0\mbox{ mod }q^2$;
$y\in Y(k)$ corresponds to $q=t^d$,
$\bar x=x^0(t)\mbox{ mod }t^{2d}$, $\bar y=y^0(t)\mbox{ mod }t^{d}$.
\hfill\qedsymbol

\noindent
{\bf Remark.} The above proof goes back to Finkelberg and Mirkovi\'c
\cite{FM}. They proved Theorem \ref{main} in the particular case that
$X$ is the affine closure of $G/N$, where $G$ is a reductive group
and $N$ is the unipotent radical of a Borel subgroup $B\subset G$.
Their main idea is to choose an opposite Borel subgroup 
$B_-\subset G$ such that the image of $\gamma_0$ has non-empty
intersection with the open orbit of $B_-$ in $X$. Our proof makes an
implicit use of a certain smooth groupoids $\Gamma_r$, $r\ge 1$, of
dimension $n=\dim X$ acting on $X$ so that there exists a closed
$\Gamma_r$-invariant subscheme $Y\subset X$ with the following 
properties: $Y\subset X$ is defined by one equation, the action of
$\Gamma_r$ on $X\setminus Y$ is transitive, and the image of
$\gamma_0$ is not contained in $Y$. More precisely, our
$X\subset {\rm Spec}\, k[x_1,\ldots,x_n,y_1,\ldots,y_l]$ is assumed
to be a complete intersection $p_1=\ldots =p_l=0$ and 
$Q:=\det\frac{\partial p}{\partial y}$. I will explicitly describe
$\Gamma_r$ in a subsequent version of this paper; here I am just
going to indicate that the corresponding Lie algebroid is generated
by the vector fields $v_i$ on $X$ such that
$v_i(x_j)=Q^r\delta_{ij}$, , $1\le i\le n$.


\newpage

\begin{thebibliography}{FM}
\bibitem[FM]{FM}
M. Finkelberg, I. Mirkovi\'c, {\it Semiinfinite flags. I. Case of global
curve $P^1$}, Differential topology, infinite-dimensional
Lie algebras, and applications, 81--112, 
Amer. Math. Soc. Transl. Ser. 2, 194, Amer. Math. Soc., Providence,
RI, 1999. See also e-print alg-geom/9707010.

\bibitem[GK]{GK}
M.~Grinberg and D.~Kazhdan, {\it Versal deformations of formal arcs},
Geometric and Functional Analysis, vol. 10 (2000), no. 3, 543-555.
See also e-print math.AG/9812104.


\end{thebibliography}
\end{document}